\documentclass[11pt]{article}

\usepackage[margin=1in]{geometry} 
\usepackage{amsmath} 
\usepackage{amssymb} 
\usepackage{graphicx} 
\usepackage{textcomp} 
\usepackage{listings} 
\usepackage{xcolor}   
\usepackage{hyperref} 
\hypersetup{
    colorlinks=true, 
    linkcolor=blue,  
    citecolor=blue,  
    filecolor=magenta, 
    urlcolor=blue   
}

\begin{document}

\title{A Technical Survey of Sparse Linear Solvers in Electronic Design Automation}

\author{
  Nityanand Rai\thanks{Work partially done while affiliated with Cadence Design Systems.} \thanks{Amazon Web Services (AWS), Seattle, WA, USA.} \\
  (ORCID: 0009-0002-2196-5023)
}

\date{September 2024} 

\maketitle

\begin{abstract}
Sparse linear system solvers ($Ax=b$) represent critical computational kernels within Electronic Design Automation (EDA), underpinning the simulation of complex physical phenomena vital for modern integrated circuit and system design. Applications such as power integrity verification and electrothermal analysis rely fundamentally on these algorithms to solve the large-scale, sparse algebraic systems derived from Modified Nodal Analysis (MNA) of circuits or Finite Element/Volume Method (FEM/FVM) discretizations of governing partial differential equations (PDEs). Problem dimensions routinely reach $10^6-10^9$ unknowns, escalating towards $10^{10}$ or more for full-chip power grid simulations \cite{Tsinghua21}, imposing stringent demands on solver scalability, memory footprint, and computational efficiency. This paper provides a technical survey of the predominant sparse solver paradigms utilized in EDA: direct methods based on matrix factorization (LU, Cholesky), iterative methods employing Krylov subspace projections (CG, GMRES, BiCGSTAB), and multilevel multigrid techniques. We examine their mathematical foundations, including convergence properties and sensitivity to matrix conditioning, implementation aspects such as sparse storage formats (CSR, CSC) and fill-in mitigation strategies (reordering algorithms), the critical role of preconditioning in accelerating iterative method convergence for ill-conditioned systems \cite{SaadIterative, ComparisonSolversArxiv}, and the potential for optimal $O(N)$ complexity offered by multigrid \cite{TrottenbergMG}. A crucial factor influencing solver choice is the performance impact of frequent matrix updates (e.g., in transient or non-linear simulations), where iterative and multigrid methods often amortize costs more effectively than direct methods requiring repeated factorization \cite{SaadIterative}. We analyze the trade-offs concerning runtime complexity, memory requirements, numerical robustness, parallel scalability (MPI, OpenMP, GPU), and numerical precision (FP32/FP64). The integration of these solvers into sophisticated EDA tools addressing system-level multiphysics challenges is discussed, illustrated by pseudocode for core algorithms. The survey concludes by emphasizing the ongoing evolution and indispensable nature of sparse solvers in enabling the design and verification of complex electronic systems.
\end{abstract}

\noindent 
\textbf{Keywords:} Sparse Matrix, Linear Solvers, Electronic Design Automation (EDA), Direct Solvers, Iterative Solvers, Krylov Subspace Methods, Multigrid Methods, Preconditioning, Power Integrity, Thermal Analysis, Computational Fluid Dynamics (CFD), Modified Nodal Analysis, Finite Element Method, Finite Volume Method

\section{INTRODUCTION}
Electronic Design Automation (EDA) software for circuit simulation (e.g., SPICE), power integrity (IR drop) analysis, thermal simulation, and electromagnetic extraction frequently necessitates the solution of large-scale sparse linear systems of equations, canonically represented as $Ax=b$. The system matrix $A \in \mathbb{R}^{N \times N}$ (or $\mathbb{C}^{N \times N}$ for AC analysis) typically originates from circuit network equations formulated via Modified Nodal Analysis (MNA) or from the spatial discretization (using Finite Element, Finite Difference, or Finite Volume Methods - FEM, FDM, FVM) of partial differential equations (PDEs) governing physical phenomena such as heat diffusion (parabolic PDE) or electrostatic potential (elliptic PDE) \cite{SaadIterative}. Modern System-on-Chip (SoC) and 3D-IC designs feature highly complex power delivery networks (PDNs) and intricate geometries, resulting in system matrices where the dimension $N$ commonly ranges from $10^6$ to $10^9$, and can escalate towards $10^{10}$ or more for full-chip, detailed PDN simulations \cite{Tsinghua21}. The computational cost associated with solving these sparse systems often dominates the overall simulation time, making the selection and optimization of the linear solver algorithm a critical design consideration, balancing trade-offs between execution speed, memory consumption (footprint), numerical stability, parallel scalability, and robustness \cite{DavisDirect, SaadIterative}.

A notable trend is the extension of EDA capabilities from chip-centric analysis towards comprehensive system-level simulation, prominently featuring electrothermal co-simulation \cite{CelsiusDatasheet, CelsiusProductPage}. This evolution is motivated by the need to accurately model the thermal behavior of entire electronic products (e.g., smartphones, laptops, data center servers), where heat generated by high-power ICs (often derived from power integrity analysis results \cite{Celsius3DICSys}) significantly affects system performance, reliability, and energy efficiency. Analyzing thermal transport through the package, PCB, enclosure, and interaction with cooling subsystems (e.g., heat sinks, fans, liquid cooling loops) necessitates coupling electrical/thermal EDA simulations with Computational Fluid Dynamics (CFD) \cite{CelsiusProductPage}. Consequently, major EDA vendors are integrating advanced CFD solvers \cite{CelsiusLaunch}, introducing the specific numerical challenges of fluid flow (e.g., solving Navier-Stokes equations) into the EDA ecosystem, often involving loosely or tightly coupled multiphysics simulations. The convergence of extreme problem scale, complex multiphysics interactions, and the need for efficient parallel execution on modern hardware makes a contemporary survey of applicable sparse solver techniques, specifically within the EDA context, particularly relevant.

High-performance EDA tools incorporate sophisticated sparse matrix solver libraries and custom implementations to manage this computational burden \cite{VoltusBlogNov13, CelsiusProductPage}. This paper presents a technical survey of the three principal classes of sparse linear solvers relevant to EDA and associated CFD domains: Direct Methods, Iterative Methods, and Multigrid Methods. We explore their underlying mathematical principles, algorithmic details, computational complexity, memory requirements, parallelization strategies (distributed and shared memory), and practical integration within modern analysis frameworks \cite{DavisDirect, SaadIterative, TrottenbergMG}, providing insights into the rationale governing their selection and deployment.

\section{FUNDAMENTALS OF SPARSE MATRICES IN EDA}

\subsection{Origin and Characteristics}
The prevalence of sparse matrices in EDA stems from the inherent locality of interactions in both electrical circuits (MNA) and discretized physical domains (FEM/FVM/FDM) \cite{DavisDirect}. The resulting system matrix $A$ typically possesses a number of non-zero elements (NNZ) that scales linearly or near-linearly with the number of unknowns $N$, i.e., $NNZ = O(N)$, in stark contrast to dense matrices where $NNZ = N^2$ \cite{DavisDirect}. The specific properties of $A$ profoundly influence solver selection and performance:
\begin{itemize}
   \item \textbf{Sparsity Pattern:} The distribution of non-zero entries, often represented by the matrix's adjacency graph, reflects the connectivity of the circuit or the discretization stencil. Patterns can range from banded (in simpler 1D discretizations) to highly irregular (in complex 3D geometries or MNA).
   \item \textbf{Condition Number ($\kappa(A)$):} Defined as $\kappa(A) = \|A\| \|A^{-1}\|$ (using an appropriate matrix norm), this number quantifies the sensitivity of the solution $x$ to perturbations in $A$ or $b$. Power grid matrices are often severely ill-conditioned, with $\kappa(A)$ reaching $10^6 - 10^{10}$ or higher \cite{ComparisonSolversArxiv}, primarily due to large variations in conductance values. High condition numbers drastically impede the convergence rate of standard iterative methods, mandating the use of powerful preconditioning techniques \cite{SaadIterative, ComparisonSolversArxiv}.
   \item \textbf{Symmetry and Definiteness:} DC power grid analysis and steady-state heat conduction problems often yield Symmetric Positive Definite (SPD) matrices, for which specialized and highly efficient solvers like Cholesky factorization or the Conjugate Gradient (CG) method exist \cite{DavisDirect}. AC analysis introduces complex arithmetic and non-Hermitian matrices. CFD momentum (Navier-Stokes) and coupled multiphysics problems frequently lead to non-symmetric, and potentially indefinite, systems \cite{SaadIterative}, requiring more general solvers like LU decomposition or GMRES/BiCGSTAB.
\end{itemize}

\subsection{Storage Formats}
Minimizing memory footprint and enabling efficient matrix operations (especially Sparse Matrix-Vector multiplication, SpMV) necessitates specialized storage schemes beyond simple dense arrays \cite{DavisDirect, EigenDoc}. Common formats include:
\begin{itemize}
   \item \textbf{Coordinate (COO):} Stores triplets \texttt{(row\_index, col\_index, value)} for each non-zero. Simple to construct but inefficient for SpMV or row/column access.
   \item \textbf{Compressed Sparse Row (CSR):} Uses three arrays: one for non-zero values, one for column indices, and one (\texttt{row\_ptr}) indicating the start of each row in the other two arrays. Optimized for fast SpMV ($y \leftarrow Ax$) and row operations.
   \item \textbf{Compressed Sparse Column (CSC):} Analogous to CSR but compressed by column. Optimized for $y \leftarrow A^T x$ and column operations. Useful for LU factorization \cite{DavisDirect}.
   \item \textbf{Others:} Variations like Block CSR (BCSR) for matrices with dense blocks, Diagonal (DIA), ELLPACK, and Hybrid formats exist, offering performance trade-offs for specific matrix structures or hardware architectures \cite{SaadIterative}.
\end{itemize}

\section{DIRECT METHODS}
Direct methods aim to compute the solution $x = A^{-1}b$ in a finite number of steps (ignoring finite-precision arithmetic effects) by explicitly or implicitly factorizing the matrix $A$ into a product of simpler (typically triangular) matrices, followed by forward and backward substitution \cite{DavisDirect}.

\subsection{LU and Cholesky Factorization}
The foundational techniques are based on Gaussian elimination:
\begin{itemize}
   \item \textbf{LU Decomposition:} For a general non-singular matrix $A$, find permutation matrices $P, Q$, a unit lower triangular matrix $L$, and an upper triangular matrix $U$ such that $PAQ = LU$ or $PA=LU$. The system $Ax=b$ is then solved via $Ly=Pb$ (forward substitution) and $UQx=y$ or $Ux=y$ (backward substitution) \cite{DavisDirect, MatlabLUDoc}. Numerical stability often requires pivoting (incorporated in $P$).
   \item \textbf{Cholesky Factorization:} If $A$ is SPD, it can be factorized as $A = LL^T$, where $L$ is a lower triangular matrix (the Cholesky factor). Optionally, permutations $P$ are used for sparsity: $PAP^T = LL^T$. This method is significantly more efficient (roughly half the computation) and numerically stable (no pivoting needed for SPD matrices) than LU decomposition \cite{DavisDirect}.
\end{itemize}
The theoretical complexity for dense matrices is $O(N^3)$ for LU and Cholesky. For sparse matrices, the complexity depends heavily on fill-in \cite{DavisDirect}.

\subsection{Key Challenge: Fill-in}
During the factorization process (Gaussian elimination steps), zero entries in the original matrix $A$ can become non-zero in the factors $L$ and $U$. This phenomenon, known as \textbf{fill-in}, can dramatically increase the number of non-zeros, leading to excessive memory consumption ($O(NNZ(L+U))$) and computational cost ($O(\text{flops}(L+U))$) \cite{DavisDirect}. For large 3D problems, fill-in can make direct methods computationally infeasible, with memory and time potentially scaling as $O(N^2)$ and $O(N^3)$ respectively in worst-case scenarios, although typically better ($O(N^{1.5})$ to $O(N^2)$) in practice for well-ordered problems \cite{DavisDirect}.

\subsection{Mitigating Fill-in: Ordering and Symbolic Analysis}
The amount of fill-in is highly sensitive to the ordering of rows and columns. Preprocessing $A$ using heuristic reordering algorithms is essential to minimize fill-in during factorization. Common strategies include \cite{DavisDirect}:
\begin{itemize}
    \item \textbf{Minimum Degree (MD):} A greedy approach that eliminates rows/columns with the fewest non-zeros at each step.
    \item \textbf{Approximate Minimum Degree (AMD):} A faster, more sophisticated variant of MD.
    \item \textbf{Nested Dissection (ND):} A graph partitioning approach (using separators) that is often asymptotically optimal for regular grids and performs well on general sparse matrices, particularly in parallel environments. Libraries like METIS or Scotch are often used for partitioning.
\end{itemize}
A \textbf{symbolic factorization} phase often precedes numerical factorization to determine the non-zero pattern of $L$ and $U$ and pre-allocate memory \cite{DavisDirect}.

\subsection{Advanced Techniques and Libraries}
Modern high-performance direct solvers employ sophisticated techniques to improve efficiency and exploit parallelism \cite{DavisDirect, ComparisonSolversArxiv, EigenDoc}:
\begin{itemize}
    \item \textbf{Supernodal/Multifrontal Methods:} These methods identify and process dense submatrices (supernodes or frontal matrices) within the sparse factorization, allowing the use of highly optimized Basic Linear Algebra Subprograms (BLAS), particularly level-3 BLAS (matrix-matrix operations), for improved cache utilization and performance \cite{DavisDirect}.
    \item \textbf{Numerical Pivoting:} Techniques like threshold partial pivoting are used in LU for numerical stability when diagonal entries are small or zero \cite{DavisDirect}.
\end{itemize}
Widely used libraries implementing these techniques include Intel MKL PARDISO, SuiteSparse (UMFPACK, CHOLMOD, KLU, SPQR), MUMPS (MUltifrontal Massively Parallel sparse direct Solver), and vendor-specific implementations \cite{EigenDoc, XyceDoc, DavisDirect}.

\subsection{Direct Solver Pseudocode (Simplified LU)}
\begin{lstlisting}[language={}, caption={Simplified LU Factorization and Solve (Illustrative, No Pivoting)}, label={lst:lu}]
// Simplified In-Place LU Factorization (No Pivoting)
function SimpleLU(A[1..n, 1..n])
 for k = 1 to n-1 do
   // Check for zero pivot (requires pivoting in practice)
   if A[k,k] == 0 error "Pivot is zero"
   // Compute multipliers for column k
   for i = k+1 to n do A[i,k] /= A[k,k] end for
   // Update trailing submatrix
   for j = k+1 to n do
     for i = k+1 to n do
       A[i,j] -= A[i,k] * A[k,j] // Rank-1 update
     end for
   end for
 end for
 // A now stores L (below diag) and U (diag and above)
end function

// Solve Ax=b using LU factors (stored in A)
function SolveWithLU(LU_factors, b) -> x
 // Forward Substitution: Solve Ly = b (unit diagonal L)
 y[1..n];
 for i = 1 to n do
   y[i] = b[i]
   for j = 1 to i-1 do y[i] -= LU_factors[i,j] * y[j] end for
 end for
 // Backward Substitution: Solve Ux = y
 x[1..n];
 for i = n down to 1 do
   x[i] = y[i]
   for j = i+1 to n do x[i] -= LU_factors[i,j] * x[j] end for
   x[i] /= LU_factors[i,i] // Divide by diagonal element of U
 end for
 return x
end function
\end{lstlisting}

\subsection{Performance and Use Cases: Intuition for Selection}
Direct solvers are favored when \textbf{maximum numerical robustness} is required, or for problems where multiple right-hand sides ($b$) need to be solved with the same matrix $A$ (factorization is done once, solves are cheap). They are the standard for traditional SPICE simulation due to their reliability, especially for \textbf{small to moderate $N$} (e.g., $N \lesssim 10^6 - 10^7$) \cite{DavisDirect}. However, the potentially high memory ($O(NNZ(L+U))$) and factorization time ($O(\text{flops}(L+U))$) cost due to fill-in limits their scalability for very large EDA problems. A critical drawback is that if the matrix $A$ changes (e.g., in non-linear iterations or transient simulations), the \textbf{entire, computationally expensive factorization phase must typically be repeated} \cite{SaadIterative}.

\section{ITERATIVE METHODS}
Iterative methods construct a sequence of approximate solutions $x^{(k)}$ that converge towards the true solution $x$, i.e., $\lim_{k \to \infty} x^{(k)} = x$. They typically require less memory than direct methods, often only $O(NNZ)$ plus storage for a few vectors ($O(N)$), making them attractive for very large systems. Their core computational kernel is usually the Sparse Matrix-Vector Multiplication (SpMV), which exhibits good parallel scalability but can be memory-bandwidth limited \cite{SaadIterative}. Convergence is achieved when a specific criterion is met, e.g., the relative residual norm $\|b - Ax^{(k)}\|_2 / \|b\|_2$ falls below a specified tolerance $\epsilon$.

\subsection{Stationary Methods}
Classic methods like Jacobi, Gauss-Seidel, and Successive Over-Relaxation (SOR) update the solution based on a simple splitting of $A$. Their convergence rate is linear and often depends on the spectral radius $\rho(G)$ of their iteration matrix $G$ (convergence requires $\rho(G) < 1$). While generally too slow for solving large EDA systems directly, they are valuable as computationally inexpensive \textbf{smoothers} within multigrid methods \cite{TrottenbergMG, SaadIterative} or as simple preconditioners.

\subsection{Krylov Subspace Methods}
These are the workhorses of modern iterative solvers. They seek an approximate solution $x^{(k)}$ within the Krylov subspace $\mathcal{K}_k(A, r_0) = \text{span}\{r_0, Ar_0, A^2r_0, \dots, A^{k-1}r_0\}$, where $r_0 = b - Ax^{(0)}$ is the initial residual. They typically exhibit superlinear convergence for well-conditioned systems. Key methods include \cite{SaadIterative}:
\begin{itemize}
   \item \textbf{Conjugate Gradient (CG):} The optimal Krylov method for SPD systems. It minimizes the A-norm of the error $\|x - x^{(k)}\|_A = \sqrt{(x - x^{(k)})^T A (x - x^{(k)})}$ over $\mathcal{K}_k(A, r_0)$. Requires one SpMV and a few vector operations per iteration.
   \item \textbf{Generalized Minimum Residual (GMRES):} Robust method for general non-symmetric systems. It minimizes the Euclidean norm of the residual $\|b - Ax^{(k)}\|_2$ over $\mathcal{K}_k(A, r_0)$. Requires storing an orthonormal basis for $\mathcal{K}_k$, leading to increasing memory and computational cost per iteration. Often \textbf{restarted} (GMRES(m)) after $m$ iterations to limit cost, potentially slowing convergence.
   \item \textbf{Bi-Conjugate Gradient Stabilized (BiCGSTAB):} An alternative for non-symmetric systems that avoids the increasing cost of GMRES by using short-term recurrences. Often converges faster than restarted GMRES but can exhibit more irregular convergence behavior. Requires two SpMV operations per iteration.
   \item \textbf{Others:} Methods like QMR (Quasi-Minimal Residual) and CGS (Conjugate Gradient Squared) offer different trade-offs.
\end{itemize}

\subsection{The Crucial Role of Preconditioning}
For the ill-conditioned matrices common in EDA ($\kappa(A) \gg 1$), the convergence of Krylov methods can be unacceptably slow or may even stagnate. \textbf{Preconditioning} is essential to accelerate convergence. The goal is to find a matrix $M$ (the preconditioner) such that $M \approx A^{-1}$ and the transformed system ($M^{-1}Ax = M^{-1}b$ or $AM^{-1}y = b, x=M^{-1}y$) is much better conditioned, i.e., $\kappa(M^{-1}A) \ll \kappa(A)$ or $\kappa(AM^{-1}) \ll \kappa(A)$, and linear systems involving $M$ are inexpensive to solve \cite{SaadIterative}. Common preconditioning strategies include:
\begin{itemize}
   \item \textbf{Simple/Diagonal (Jacobi):} $M = \text{diag}(A)$. Very cheap but often ineffective for challenging problems.
   \item \textbf{Incomplete Factorizations (ILU/ICC):} Approximate LU (ILU) or Cholesky (ICC) factorizations where fill-in is restricted based on level (ILU(k)) or magnitude threshold (ILUT). Offer a trade-off between effectiveness and setup/application cost \cite{SaadIterative}.
   \item \textbf{Algebraic Multigrid (AMG):} A powerful class of preconditioners, particularly effective for elliptic PDE discretizations. Constructs a hierarchy of coarser representations of the matrix algebraically, without requiring geometric grid information. Can provide optimal $O(N)$ preconditioning cost \cite{TrottenbergMG}. See Section \ref{sec:multigrid}.
   \item \textbf{Domain Decomposition (DD):} Partition the problem domain/graph and solve smaller problems on subdomains, combining solutions iteratively (e.g., additive/multiplicative Schwarz). Often involves solving interface problems (Schur complement methods). Highly parallelizable \cite{SaadIterative}.
   \item \textbf{Physics/Graph-Based:} Exploit problem structure, e.g., block Jacobi based on physical components, multi-level aggregation specific to power grid topology \cite{EigenDoc, Garyfallou16}.
\end{itemize}

\subsection{Iterative Solver Pseudocode (Conjugate Gradient)}
\begin{lstlisting}[language={}, caption={Preconditioned Conjugate Gradient (PCG) Algorithm}, label={lst:cg}]
// Preconditioned Conjugate Gradient for solving Ax=b, where A is SPD
// M is the preconditioner (approximates A^-1)
function PCG(A, b, x0, M, tol) -> x
 r = b - A*x0          // Initial residual
 z = solve(M, r)       // Apply preconditioner: z = M^-1 * r
 p = z                 // Initial search direction
 rho_old = dot(r, z)
 for k = 1 to max_iter do
   if sqrt(dot(r,r)) < tol * sqrt(dot(b,b)) break // Check residual norm
   Ap = A*p            // Matrix-vector product
   alpha = rho_old / dot(p, Ap) // Step length
   x = x + alpha * p   // Update solution
   r = r - alpha * Ap  // Update residual
   z = solve(M, r)     // Apply preconditioner
   rho_new = dot(r, z) // New dot product
   if rho_new == 0 break // Avoid division by zero (exact solution found)
   beta = rho_new / rho_old // Improvement factor
   p = z + beta * p    // Update search direction
   rho_old = rho_new
 end for
 return x
end function
\end{lstlisting}

\subsection{Performance and Use Cases: Intuition for Selection}
Iterative methods are generally the preferred choice for \textbf{extremely large N} ($N > 10^7 - 10^8$) or when \textbf{memory constraints} preclude direct methods. Their performance per iteration is dominated by SpMV ($O(NNZ)$) and preconditioner application ($O(\text{cost}(M^{-1}))$). They are particularly advantageous when the matrix $A$ \textbf{changes frequently}, as the cost per solve (after potentially expensive initial preconditioner setup) can be much lower than re-factorizing \cite{SaadIterative}. Their inherent parallelism, especially in SpMV and vector operations, makes them well-suited for \textbf{massively parallel architectures} (multi-core CPUs, GPUs, clusters) \cite{Garyfallou16}. However, their \textbf{convergence behavior is problem-dependent} and critically relies on the \textbf{effectiveness of the preconditioner} \cite{SaadIterative}. Poor preconditioning can lead to slow convergence or stagnation. Mixed-precision techniques (e.g., using FP32 for SpMV and FP64 for residuals) can offer speedups \cite{OktayCarsonMixedGMRES}.

\section{MULTIGRID METHODS}\label{sec:multigrid}
Multigrid methods are a class of algorithms originally developed for solving elliptic PDEs discretized on structured grids, often exhibiting asymptotically optimal or near-optimal computational complexity, potentially $O(N)$ \cite{TrottenbergMG}. They can be used either as standalone solvers or as highly effective preconditioners for Krylov methods (e.g., AMG-PCG).

\subsection{Core Idea: Solving on Multiple Scales}
The fundamental principle of multigrid is to accelerate convergence by utilizing a hierarchy of grids (or matrix representations) to efficiently eliminate error components across different spatial frequencies. A typical multigrid cycle involves \cite{TrottenbergMG}:
\begin{enumerate}
    \item \textbf{Pre-Smoothing:} Apply a few iterations of a simple iterative method (e.g., Gauss-Seidel, Jacobi - the smoother) to the current approximation $u_h$ on the fine grid $h$. This primarily damps high-frequency error components.
    \item \textbf{Residual Computation:} Calculate the fine-grid residual $r_h = b_h - A_h u_h$.
    \item \textbf{Restriction:} Transfer the residual $r_h$ to the next coarser grid $H$, obtaining $r_H = R r_h$, where $R$ is the restriction operator (e.g., injection, full-weighting).
    \item \textbf{Coarse Grid Solve:} Solve the residual equation $A_H e_H = r_H$ on the coarse grid $H$ to find the coarse-grid error correction $e_H$. Here $A_H$ is the coarse-grid operator. This step is performed recursively, applying multigrid cycles until the coarsest grid is reached, where the system is solved directly (or with many iterations).
    \item \textbf{Prolongation/Interpolation:} Interpolate the error correction $e_H$ back to the fine grid $h$, obtaining $e_h = P e_H$, where $P$ is the prolongation operator (often related to $R$, e.g., $P=c R^T$).
    \item \textbf{Correction:} Update the fine-grid approximation: $u_h \leftarrow u_h + e_h$.
    \item \textbf{Post-Smoothing:} Apply a few more smoothing iterations to damp any high-frequency errors introduced by the prolongation step.
\end{enumerate}
Different cycle structures exist (V-cycle, W-cycle, F-cycle, Full Multigrid - FMG) based on how the recursive coarse-grid solves are organized \cite{TrottenbergMG}.

\subsection{Geometric vs. Algebraic Multigrid (AMG)}
\begin{itemize}
   \item \textbf{Geometric Multigrid (GMG):} Requires an explicit hierarchy of nested geometric grids. Restriction and prolongation operators are defined based on grid geometry (e.g., linear interpolation). Relatively intuitive to implement for \textbf{structured or semi-structured meshes} \cite{TrottenbergMG}.
   \item \textbf{Algebraic Multigrid (AMG):} Constructs the entire multigrid hierarchy (coarse "grids", restriction $R$, prolongation $P$, and coarse operators $A_H$) purely from the entries of the matrix $A_h$. Typically uses heuristics based on "strength of connection" between matrix variables to define coarse grids and interpolation weights. The coarse operator is often formed via the Galerkin projection: $A_H = R A_h P$. AMG is essential for problems defined on \textbf{complex, unstructured meshes} prevalent in EDA and CFD, where geometric grid hierarchies are unavailable or difficult to construct \cite{TrottenbergMG}. The setup phase (building the hierarchy) can be computationally intensive but is amortized over solves or preconditioner applications.
\end{itemize}

\subsection{Multigrid Pseudocode (Recursive V-Cycle)}
\begin{lstlisting}[language={}, caption={Recursive Multigrid V-Cycle Function (as Solver)}, label={lst:mg}]
// Recursive Multigrid V-Cycle to solve A_h * u_h = b_h
// h denotes the current level, H the next coarser level
// nu1, nu2 are the number of pre/post-smoothing steps
function V_Cycle(A_h, b_h, u_h, level h) -> u_h
 if h is the coarsest level then
   Solve A_h * u_h = b_h (e.g., direct solver or many iterations)
   return u_h
 else
   // 1. Pre-Smoothing
   u_h = Smooth(A_h, b_h, u_h, nu1)
   // 2. Residual Computation
   r_h = b_h - A_h * u_h
   // 3. Restriction
   r_H = Restrict(r_h, h) // R operator
   // 4. Coarse Grid Solve (Recursive Call)
   e_H = 0 // Initial guess for error on coarse grid
   A_H = GetCoarseGridOperator(A_h, h) // Often R*A_h*P
   e_H = V_Cycle(A_H, r_H, e_H, level H)
   // 5. Prolongation
   e_h = Prolongate(e_H, h) // P operator
   // 6. Correction
   u_h = u_h + e_h
   // 7. Post-Smoothing
   u_h = Smooth(A_h, b_h, u_h, nu2)
   return u_h
 end if
end function
\end{lstlisting}

\subsection{Performance and Use Cases: Intuition for Selection}
For problems where it is applicable (typically discretized elliptic PDEs like those arising in Poisson solves for pressure in CFD or DC power grids), Multigrid offers the potential for \textbf{optimal $O(N)$ or near-optimal $O(N \log N)$ computational complexity}, meaning the solution time scales linearly with problem size \cite{TrottenbergMG}. It can handle \textbf{extremely large systems} ($N > 10^9$) efficiently, often exhibiting convergence rates independent of the mesh size $h$. When $A$ \textbf{changes frequently}, the fast convergence per V-cycle makes it efficient, provided the setup cost (especially for AMG) is amortized or updates are cheap \cite{TrottenbergMG}. However, implementation, particularly of robust AMG, is significantly more \textbf{complex} than standard Krylov methods. Performance can be sensitive to the choice of smoother, restriction/prolongation operators, cycle type, and the quality of the coarse-grid operators. Specialized libraries (e.g., HYPRE, ML (Trilinos), PETSc's MG components) are commonly employed \cite{SaadIterative}.

\section{COMPARATIVE ANALYSIS: CHOOSING THE RIGHT SOLVER}
Selecting the most suitable sparse linear solver for a given EDA application involves a careful assessment of multiple factors: matrix properties (symmetry, definiteness, conditioning), problem scale ($N$, NNZ), required accuracy ($\epsilon$), available computational resources (memory, CPU/GPU cycles), parallel computing environment, and the frequency of matrix updates.

\begin{table}[h!]
  \centering 
  \caption{Comparative Overview of Sparse Solver Paradigms in EDA Context}
  \label{tab:solver_comparison}
   \small 
   \sloppy 
  \begin{tabular}{|p{0.15\textwidth}|p{0.25\textwidth}|p{0.25\textwidth}|p{0.25\textwidth}|} 
    \hline 
    \textbf{Feature} & \textbf{Direct \cite{DavisDirect}} & \textbf{Iterative (Krylov) \cite{SaadIterative}} & \textbf{Multigrid (AMG) \cite{TrottenbergMG}} \\
    \hline 
    \textbf{Best For} & Small/Med N; Max Robustness; Multiple RHS & Very Large N; Memory Constrained; Parallel Systems & Elliptic PDE-like; Max Scalability ($O(N)$); Preconditioning \\ \hline
    \textbf{Strength} & Guaranteed Sol. (if stable); Predictable Cost; Robust & Low Mem ($O(NNZ)$); Good Parallelism (SpMV); Adaptable & Optimal Time ($O(N)$ possible); Fast Convergence (h-independent) \\ \hline
    \textbf{Weakness} & High Mem/Time ($\propto NNZ(\text{Factors})$); Poor Scale ($N^{1.5}\!-\!N^2$); Fill-in & Needs Eff. Precond.; Convergence Varies; Potential Stagnation & High Complexity; Setup Cost (AMG); Problem Dependent \\ \hline
    \textbf{Setup Cost} & High (Factorization, Ordering) & Low (vectors) to High (complex PC setup) & Moderate to Very High (AMG hierarchy) \\ \hline
    \textbf{Solve Cost (per RHS)} & Low (Subst: $\propto NNZ(\text{Factors})$) & Moderate ($k \times (\text{SpMV} + \text{PC\_apply})$) & Low (few cycles, potentially $O(N)$) \\ \hline 
    \textbf{Dynamic $A$} & Very Expensive (Re-factor) & Often Efficient (Re-solve, may update PC) & Efficient (Re-solve, may update hierarchy) \\
    \hline 
  \end{tabular}
  \fussy 
\end{table}

\textbf{High-Level Heuristic:} Use \textbf{Direct} methods for guaranteed solutions on smaller problems ($N \lesssim 10^7$) or when robustness is critical and factorization cost is acceptable \cite{DavisDirect}. Use \textbf{Iterative} methods for very large N or when memory is limited, but invest heavily in finding/tuning an effective \textbf{Preconditioner} \cite{SaadIterative}. Use \textbf{Multigrid} (often AMG as a solver or preconditioner) for problems arising from elliptic PDEs demanding maximum scalability and speed, accepting higher implementation or library dependency complexity \cite{TrottenbergMG}.

\section{IMPLEMENTATION CONSIDERATIONS AND LIBRARIES}
Developing or deploying high-performance sparse solvers requires attention to numerous practical details: choosing appropriate storage formats \cite{DavisDirect, EigenDoc}, implementing robust fill-reducing orderings \cite{DavisDirect}, designing algorithms for parallel execution (MPI for distributed memory, OpenMP/CUDA/SYCL for shared memory/accelerators) targeting load balancing and minimizing communication \cite{SaadIterative, Garyfallou16, CookExascalePower}, ensuring numerical stability through pivoting or careful preconditioner design \cite{DavisDirect, SaadIterative}, managing numerical precision (FP32/FP64 trade-offs, mixed-precision solvers \cite{OktayCarsonMixedGMRES}), and selecting appropriate convergence criteria. Memory bandwidth often becomes a bottleneck for iterative methods dominated by SpMV \cite{SaadIterative}. Load balancing and minimizing communication can be particularly challenging for the highly irregular matrix structures often encountered in EDA \cite{SaadIterative}.

Leveraging highly optimized numerical libraries is standard practice in EDA tool development. Prominent examples include:
\begin{itemize}
    \item \textbf{Intel Math Kernel Library (MKL):} Includes PARDISO (direct) and iterative solvers (CG, GMRES) with preconditioners.
    \item \textbf{SuiteSparse:} A collection by T. Davis including UMFPACK (multifrontal LU), CHOLMOD (Cholesky), KLU (specialized for circuits), SPQR (sparse QR) \cite{DavisDirect}.
    \item \textbf{PETSc \& Trilinos:} Comprehensive frameworks for parallel scientific computing, offering a wide array of linear/non-linear solvers, preconditioners (including AMG via HYPRE or ML), and data structures.
    \item \textbf{MUMPS:} Parallel multifrontal direct solver.
    \item \textbf{HYPRE:} Scalable linear solvers library, strong focus on parallel AMG preconditioners.
    \item \textbf{GPU Libraries:} cuSPARSE (NVIDIA), rocSPARSE (AMD) for sparse BLAS on GPUs.
    \item \textbf{Proprietary/In-house}: EDA vendors often develop highly tuned proprietary solvers optimized for specific application characteristics (e.g., \cite{XyceDoc}).
\end{itemize}

\section{SOLVER INTEGRATION IN ADVANCED EDA AND SYSTEM SIMULATION TOOLS}
Modern EDA tools deploy these solver techniques in sophisticated ways, often using hybrid approaches tailored to specific analysis needs.

\subsection{Power Integrity Analysis}
Simulating multi-billion node PDNs for DC (IR drop) or transient analysis requires extreme scalability \cite{Tsinghua21, CookExascalePower}.
\begin{itemize}
    \item \textbf{DC Analysis:} Often results in large, SPD (or nearly SPD) systems. Massively parallel \textbf{AMG-preconditioned CG (AMG-PCG)} is a common state-of-the-art approach, leveraging multigrid's scalability as a preconditioner \cite{CookExascalePower}. Hierarchical or domain decomposition methods are also employed to manage complexity \cite{Tsinghua21}.
    \item \textbf{Transient Analysis:} Involves time-stepping, potentially changing matrix entries if non-linear devices are present. Iterative methods (PCG, BiCGSTAB) are often favored due to the cost of repeated direct factorization, possibly with preconditioners that can be efficiently updated \cite{SaadIterative}.
\end{itemize}
Massive parallelism (thousands of cores) is essential \cite{Tsinghua21, CookExascalePower}.

\subsection{Thermal and Multiphysics Analysis (including CFD)}
System-level electrothermal tools couple electrical models with heat transfer (solids) and potentially fluid dynamics (CFD) for cooling \cite{CelsiusProductPage, CelsiusDatasheet, Celsius3DICSys}. This requires solving different types of equations, often demanding different solvers within a coupled iteration framework:
\begin{itemize}
   \item \textbf{Pressure Equation (CFD):} Arises from enforcing mass conservation (e.g., in SIMPLE/PISO pressure-correction schemes for incompressible flow), resulting in a large, sparse, often symmetric (or nearly symmetric) Poisson-like equation for pressure correction. \textbf{AMG} is typically the most efficient solver or preconditioner for this elliptic problem due to its scalability \cite{TrottenbergMG}.
   \item \textbf{Momentum Equations (CFD):} Discretized Navier-Stokes equations lead to large, non-symmetric, convection-dominated systems. Stabilized \textbf{Iterative Krylov methods (GMRES, BiCGSTAB)} are commonly used, often preconditioned with relatively simpler methods like \textbf{ILU(0), Block Jacobi}, or sometimes domain decomposition, due to the complexity and changing nature of the flow field \cite{SaadIterative}.
   \item \textbf{Energy/Heat Equation (Thermal/CFD):} A scalar convection-diffusion equation. Depending on whether it's diffusion or convection dominated, \textbf{PCG (if SPD), BiCGSTAB, or GMRES} are used, potentially preconditioned with ILU, or even \textbf{AMG} if diffusion dominates or high scalability is needed \cite{SaadIterative, TrottenbergMG}.
\end{itemize}
These solvers operate within parallel CFD frameworks, iteratively exchanging boundary conditions or source terms between the electrical, thermal conduction (solid), and fluid dynamics (coolant) domains until convergence \cite{Celsius3DICSys}. The integration of CFD significantly broadens the range of sparse solver requirements within the EDA landscape.

\section{CONCLUSION}
Sparse linear system solvers form the computational backbone of numerous critical EDA and system-level simulation tasks. Direct methods, while robust, face inherent scalability limitations due to factorization cost and fill-in, making them suitable primarily for smaller problems or specific contexts like traditional SPICE \cite{DavisDirect}. Iterative methods offer superior scalability in terms of memory and parallelism, particularly advantageous for dynamic problems where matrices change, but their performance critically hinges on the availability and effectiveness of preconditioning strategies to handle ill-conditioned systems common in EDA \cite{SaadIterative, ComparisonSolversArxiv}. Multigrid methods, especially AMG, provide the potential for optimal $O(N)$ complexity for elliptic-dominated problems (e.g., DC power grids, CFD pressure solves), delivering maximum scalability, albeit at the cost of increased implementation complexity \cite{TrottenbergMG}.

The selection of an appropriate solver necessitates a nuanced understanding of the trade-offs between computational complexity (factorization vs. iteration cost, setup vs. solve cost), memory footprint, numerical robustness, sensitivity to matrix properties (conditioning, symmetry), parallel scaling characteristics, and the impact of matrix dynamism. In multiphysics simulations, hybrid approaches employing different specialized solvers for distinct physical equations are standard. State-of-the-art EDA tools achieve performance breakthroughs by integrating massively parallel implementations of these diverse solver paradigms \cite{Garyfallou16, Tsinghua21, CookExascalePower}. Continued research and development in sparse solver algorithms, preconditioning techniques, parallel implementation strategies, adaptation to emerging hardware architectures, and potentially AI/ML-guided solver selection remain vital for enabling the design and verification of increasingly complex next-generation electronic systems.


\end{document}